\documentclass{svmult}
\usepackage[utf8]{inputenc}

\usepackage{makeidx}         
\usepackage{graphicx}        
\usepackage{multicol}        
\usepackage[bottom]{footmisc}

\makeindex             

\usepackage{float}
\usepackage{amsfonts}
\usepackage{amsmath}
\usepackage{mathtools}
\usepackage{hyperref}
\usepackage[hypcap=true,list=true]{subcaption}
\hypersetup{
    unicode=false,          
    pdftitle={Efficient Mesh Optimization Using the Gradient Flow of the Mean Volume},    
    pdfkeywords={mesh} {polyhedron} {smoothing} {regularization}, 
    colorlinks=true,       
    pdfborder=0 0 0,
    linkcolor=blue,          
    citecolor=blue,        
    filecolor=blue,      
    urlcolor=blue           
}

\usepackage{rotating}

\usepackage{xcolor}
\usepackage{lineno}

\usepackage{booktabs} 
\usepackage{array} 
\usepackage{paralist} 
\setdefaultleftmargin{1.6em}{2.2em}{1.87em}{1.7em}{1em}{1em}
\usepackage{verbatim} 

\usepackage{fancyhdr} 
\usepackage{todonotes}

\newcommand{\R}{\mathbb{R}}

\newcommand{\iq}{\operatorname{iq}}
\newcommand{\mr}{\operatorname{mr}}

\newcommand{\vol}{\operatorname{vol}}
\newcommand{\area}{\operatorname{area}}
\newcommand{\perim}{\operatorname{perim}}

\newcommand{\edgelength}{\operatorname{perim}}
\def\co{\colon\thinspace}

\begin{document}

\title*{Laplacian smoothing revisited}

\author{Dimitris Vartziotis\inst{1,2,3} \texorpdfstring{and}{\and} Benjamin Himpel\inst{1}}

\institute{
  \it TWT GmbH Science \& Innovation, Department for Mathematical Research \& Services,
  Bernh{\"a}user Stra{\ss}e~40--42,
  73765~Neuhausen, Germany \and
  NIKI Ltd.\ Digital Engineering, Research Center,
 205~Ethnikis Antistasis Street, 45500~Katsika, Ioannina, Greece \and
 Corresponding author. E-mail address: dimitris.vartziotis@nikitec.gr
}




\maketitle

\begin{abstract}
We show that the driving force behind the regularizing effect of Laplacian smoothing on surface elements is the popular mean ratio quality measure. We use these insights to provide natural generalizations to polygons and polyhedra. The corresponding functions measuring the quality of meshes are easily seen to be convex and can be used for global optimization-based untangling and smoothing. Using a simple backtracking line-search we compare several smoothing methods with respect to the resulting mesh quality. We also discuss their effectiveness in combination with topology modification.
\end{abstract}

\section{Introduction}

Mesh smoothing is the process of changing vertex positions in a mesh in order to improve the mesh quality for finite element analysis \cite{StrangFix2008}. Mesh smoothing methods can be classified \cite{MeiTipper2013_SmoothingPlanarMeshes,Owen1998,Wilson2011_HexahedralSmoothing} as geometry-based \cite{Field1988,VartziotisWipperSchwald2009}, optimization-based \cite{Branets2005, FreitagJonesPlassman1995,Parthasarathy1991,Shivanna2010,LengZhangXu_NovelGeometricFlowDrivenApproach2012,SastryShontz2012_MeshQualityImprovement,BrewerDiachinKnuppLeurentMelander2003}, phyics-based \cite{Shimada2000_Bubbles} and combinations thereof \cite{CanannTristanoStaten1998,Freitag1997,ChenTristanoKwok2003}. Effective mesh improvement procedures used in industrial applications generally require topology modifications of the mesh by adding or removing vertices and changing the connectivity of the vertices given by the edges while maintaining the mesh geometry \cite{BossenHeckbert1996,FreitagOllivierGooch1997,
KlingnerShewchuk2007}. Sometimes an untangling preprocessing step is necessary, which changes vertex positions in a mesh in order to get a valid mesh consisting of non-inverted elements. Untangling methods often coincide with optimization-based smoothing methods \cite{Knupp2001EWC,Li2000,FreitagPlassmann2000,AmentaBernEppstein1999}. To date there does not exist an algorithm which guarantees an untangled mesh, even if it is known to have an untangling solution.

An optimization-based smoothing method tries to maximize a mesh quality measure. Such a mesh quality usually consists of an element quality measure and a way to combine these to a global measure on the mesh. A very popular class of measures are algebraic quality measures \cite{Knupp2001}, and in particular the mean ratio quality measure. Some of these algebraic quality measures are used for the global optimization-based smoothing methods implemented in the Mesh Quality Improvement Toolkit Mesquite \cite{BrewerDiachinKnuppLeurentMelander2003}. A measure mostly found in the mathematical literature is the isoperimetric quotient, which can be related to the mean ratio quality measure.

Due to its simplicity and speed, Laplacian smoothing is by far the most popular smoothing method. It can be derived from a finite difference approximation of the Laplace operator \cite{BuellBush1973}. Even though it is the most run-time efficient method, it is not considered an effective method with respect to any popular quality measure. The class of geometric element transformation methods (GETMe) were introduced in \cite{VartziotisAthanasiadisGoudasWipper2008} in order to bridge the gap between efficient and effective methods. Extensions from surfaces to finite element volume meshes and generalizations to other polygon and volume types have been analyzed both theoretically and numerically in a series of papers. See \cite{VartziotisWipperPapadrakakis2013,VartziotisHimpel2014} and the references therein. In particular, we have seen in \cite{VartziotisHimpel2014} that the tetrahedral GETMe smoothing introduced in \cite{VartziotisWipperSchwald2009} can be viewed as the gradient ascent of the volume function 
and can be used for untangling.

It turns out that Laplacian smoothing on surface meshes maximizes a concave quality function, which is intimately related to the mean ratio quality measure. This concave function can be generalized to polygons and polyhedra. The gradient ascent of these measures yield efficient untangling and smoothing methods, which lend themselves to incorporating topology modifications, because worse elements tend to be smaller in the resulting mesh and can possibly be removed more easily while preserving the geometry. 

After some preliminaries on gradients of geometric functions in Section 2 we present the isoperimetric quotient and a derived concave mesh quality function in Section 3. In Section 4 we review Laplacian smoothing, relate it to a mesh quality function and suggest generalizations. Section 5 discusses ways of enhancing the smoothing methods by weights and topology modifications.

\section{Gradient vector fields}

In order to systematically construct GETMe smoothing methods based on gradient ascent and in order to relate it to Laplacian smoothing, it is essential to give closed formulas for the gradient of some basic geometric functions. In this section we introduce the notation and present the formulas. 

\subsection{Basic geometric functions}

For a function $f\co \R^n \mapsto \R$ the gradient vector field is a map $\R^n \to \R^n$ given by
\[
 \nabla f = \left(\frac{\partial}{\partial x_1} f,\ldots,\frac{\partial}{\partial x_n} f\right) \;.
\]
This yields the gradient ascent transformation
\[
 x \mapsto x + \nabla f_x \quad \text{for } x\in \R^n\;.
\]
Depending on the situation, the transformation can be parametrized by a scaling parameter $\rho$
\[
 x \mapsto x + \rho\nabla f_x \quad \text{for } x\in \R^n.
\]

Let us denote by
\[
 x = (x_1,\ldots,x_n) \in \R^{3n} \quad \text{for }x_i \in \R^3
\]
the geometric element determined by its vertices $(x_1,\ldots,x_n)$. This could be a tetrahedron given by 4 vertices, a triangle given by 3 vertices, more complicated elements and even entire meshes determined by $n$ vertices together with their connectivity. Simple geometric measures like
\begin{align*}
& \text{the volume of a tetrahedron} \quad \vol(x) = \frac{1}{6}((x_2-x_1)\times (x_3-x_1))\cdot (x_4-x_1),\\
& \text{the surface area of a triangle} \quad \area(x) = \frac{1}{2} \|(x_2 - x_1) \times (x_3 - x_1)\|, \text{ and}\\
& \text{the perimeter of a polygon} \quad \edgelength(x) = \sum_{i=1}^n \|x_{i+1}-x_i\|, 
\end{align*}
as well as certain combinations, variations and generalizations thereof, yield simple geometric transformations for polygons and polyhedra by considering their gradients. Positive volumes correspond to positively oriented tetrahedra. If we are in the Euclidean plane we will use the determinant in order to have a signed area function. 

Consider the vector $\nu_x/\|\nu_x\|$ orthonormal to the oriented triangle $x=(x_1,x_2,x_3)$ where
\begin{equation*}
\nu_x = (x_2 - x_1) \times (x_3 - x_1)\;.
\end{equation*}
If $x=(x_1,\ldots,x_n)$ let us write
\[
 \nu(i,j,k) = \nu_{(x_i,x_j,x_k)}\;.
\]
Then a straight-forward computation yields the formulas in Table \ref{table_formulas} for gradients of the above simple geometric functions.
\begin{table}[ht]
\centering
\begin{tabular}{| >{$\displaystyle}l<{$} | >{$\displaystyle}l<{$} |}
\hline
\centering
\rule{0pt}{2.7ex}f(x) & \nabla f\\[0.5ex]
\hline
\rule{0pt}{4ex}
\vol(x) & \frac{1}{6}(\nu_(4,3,2),\nu(4,1,3),\nu(4,2,1),\nu(1,2,3))\\
\rule{0pt}{4.5ex}
\area(x) &  \frac{1}{2}\left((x_2-x_3)\times \frac{\nu_x}{\|\nu_x\|},(x_3-x_1)\times \frac{\nu_x}{\|\nu_x\|},(x_1-x_2)\times \frac{\nu_x}{\|\nu_x\|} \right)\\
\rule{0pt}{4.5ex}
\edgelength(x) & \left(\frac{x_1-x_n}{\|x_1-x_n\|}+\frac{x_1-x_2}{\|x_1-x_2\|},\ldots,\frac{x_n-x_{n-1}}{\|x_n-x_{n-1}\|}+\frac{x_n-x_1}{\|x_n-x_1\|}\right)\\
\rule{0pt}{3.5ex}
\|x\|^n & n x \|x\|^{n-2}\\[0.5ex]
\hline
\end{tabular}
\caption{The gradients of some important geometric functions\label{table_formulas}}
\end{table}

\subsection{Generalizations}

This section is a bit technical and can be skipped on the first read, but has been provided for completeness. The above geometric functions for triangles and tetrahedra generalize to polygons and polyhedra. If $(x_1,\ldots,x_{n})$ are the vertices of a closed polygonal curve in $\R^3$, then
\begin{align*}
 \nu(1,\ldots,n) &\coloneqq \sum_{j=2,\ldots,n} (x_j - x_1) \times (x_{j+1} - x_1)\\
 & =  x_1 \times x_2 + x_2 \times x_3 + \ldots + x_{n-1} \times x_n + x_n \times x_1
\end{align*}
is the normal vector of the oriented curve. This allows us to define the area enclosed by the polygonal curve in $\R^3$
\[
\area(x) = \frac{1}{2} \|\nu(1,\ldots,n)\|.
\]
This function on surface elements extends linearly to entire surface meshes in two dimensions and three dimensions. If a vertex $x_0$ on a planar surface mesh is an inner node, then $\frac{\partial}{\partial x_0} \area$ vanishes. On a feature surface $\frac{\partial}{\partial x_0} \area$ will not move much. If $x_0$ lies on a boundary edge
\begin{equation*}
\frac{\partial}{\partial x_0} \area  = \frac{1}{2}(x_{1}-x_{n})\times \frac{\nu(0,1,n)}{\|\nu(0,1,n)\|}\;,\\
\end{equation*}
where $(x_1,\ldots,x_n)$ is the curve corresponding to the link of $x_0$. Notice, that we assume, that $x_0$ lies on a single boundary curve.

If $x_0$ is an inner node and the only free variable of the mesh, then $\nabla \vol$ vanishes. While for inner vertices, that is, vertices which lie in the interior of the surface, the definition of a volume function is of no consequence, it is a subtle issue when $x_0$ is a boundary node. If the boundary surface consists of only triangles and we consider the link  of $x_0$ given by a closed curve corresponding to the vertices $(x_1,\ldots,x_n)$, then
\begin{equation*}
\nabla\vol(x) = \frac{1}{6}\nu(1,\ldots,n)\;, \text{ when $x_0$ is the only free variable.}
\end{equation*}
If the surface consists of polygons rather than just triangles, we need to consider different triangulations of these polygons. Given a regular polygon with $n$ vertices, the number of triangulations is given by
\[C_k = \frac{1}{k+1} \begin{pmatrix}2k\\k \end{pmatrix}\;,\]
where $k = n-2$. Each face will then contribute an average of normal vectors associated to the sum of adjacent triangles. For example, the pentagon (in the case of a dodecahedron) has 5 triangulations. If $x_0$ is the only free vertex with the boundary of the adjacent three pentagons corresponding to the closed curve $(x_1,\ldots,x_9)$ with $x_1$, $x_4$ and $x_7$ being connected to $x_0$ via an edge, then 
\[
 \nabla\vol = \frac{1}{6}\cdot\frac{1}{5} (2\nu(1,4,7) + \nu(1,\ldots,9) + \nu(1,2,4,5,7,8) + \nu(1,3,4,6,7,9)).
\]
We get particularly nice formulas, if there is some symmetry. However, for computational purposes, this symmetry is not necessary. We just need to consider the contributions of all the faces. Let $f_n = (x_1,\ldots,x_n)$ be a face with $n$ with $x_1$ being the free node. Then the contribution $c_n$ by $f_n$ to $\nabla \vol$ is 
\begin{align*}
 c_3 = & \nu(1,2,3)\\
 c_4 = & \frac{1}{2}(\nu(1,2,4) + \nu(1,2,3,4))\\
 c_5 = & \frac{1}{5}(2\nu(1,2,5) + \nu(1,2,3,5) + \nu(1,2,4,5) + \nu(1,2,3,4,5)) \\
 c_6 = & \frac{1}{6}(5\nu(1,2,6) + 2\nu(1,2,3,6) + \nu(1,2,3,4,6)\\
 &{} +\nu(1,2,3,5,6) + \nu(1,2,3,4,5)) \end{align*}
There is a recursive formula, but as we will not be using this in this paper, we will continue with more interesting issues.

\section{The isoperimetric quotient}

The classical isoperimetric problem is to determine a planar figure of the largest possible area whose perimeter has a specific length. The proximity of a shape to the solution of the isoperimetric problem is measured by the two-dimensional isoperimetric quotient
\[
 \iq_2(x) = \frac{\area(x)}{C\perim(x)^2},
\]
where $C$ is chosen so that $\max(\iq_2) = 1$. If we only allow polygons with a fixed number of vertices, then the regular polygons maximize $\iq_2$. 
This problem becomes harder in higher dimensions. In fact, there does not exist a mathematical proof, that the regular icosahedron maximizes the three-dimensional isoperimetric quotient for icosahedra \cite{bezdek2007}. Since the volume comes with a useful sign, we will denote by $\iq_3$ the square root of the usual isoperimetric quotient
\[\iq_3(x) = \frac{\vol(x)}{C\area(x)^{3/2}},\]
where $\area$ is the area of the boundary surface of the polyhedron. Again, this quotient measures how round these objects are.

\subsection{Two-dimensional global optimization-based smoothing}

\begin{figure}[ht]
\centering
\def\svgwidth{\textwidth/3-1.5mm}
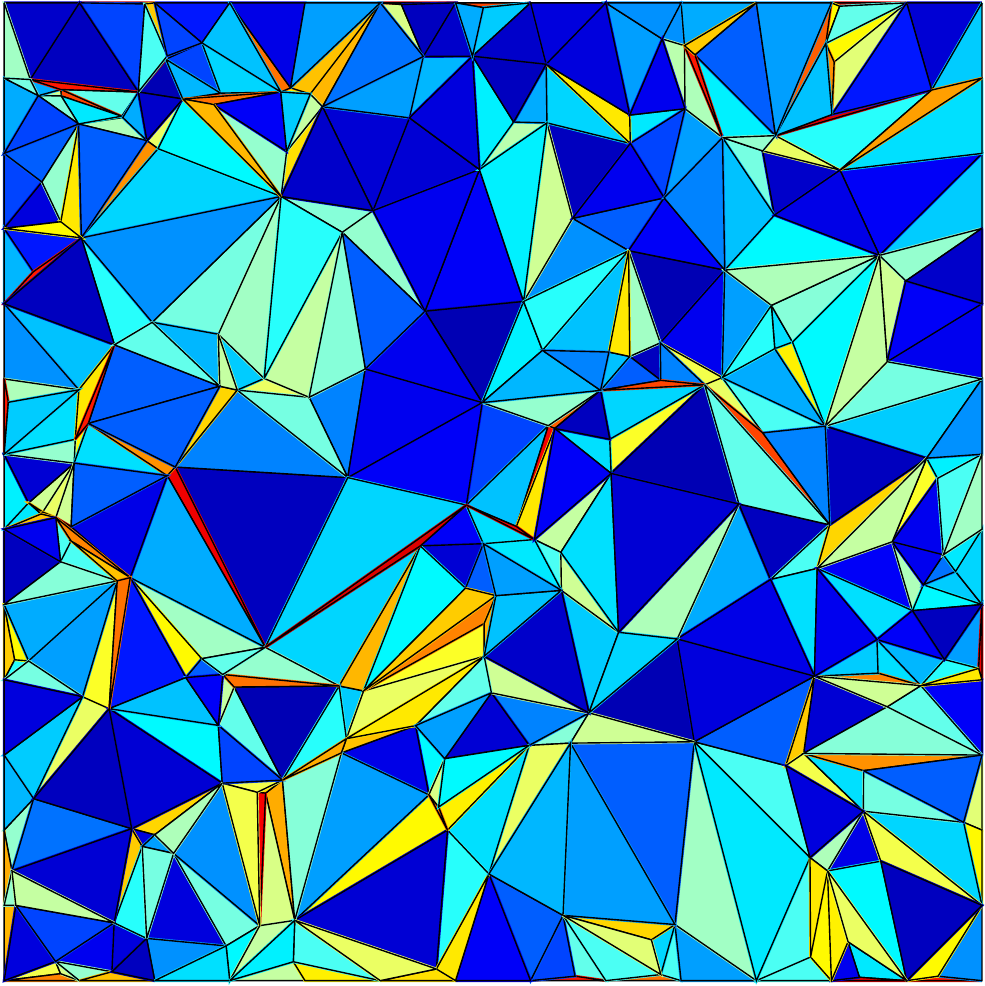
\def\svgwidth{\textwidth/3-1.5mm}
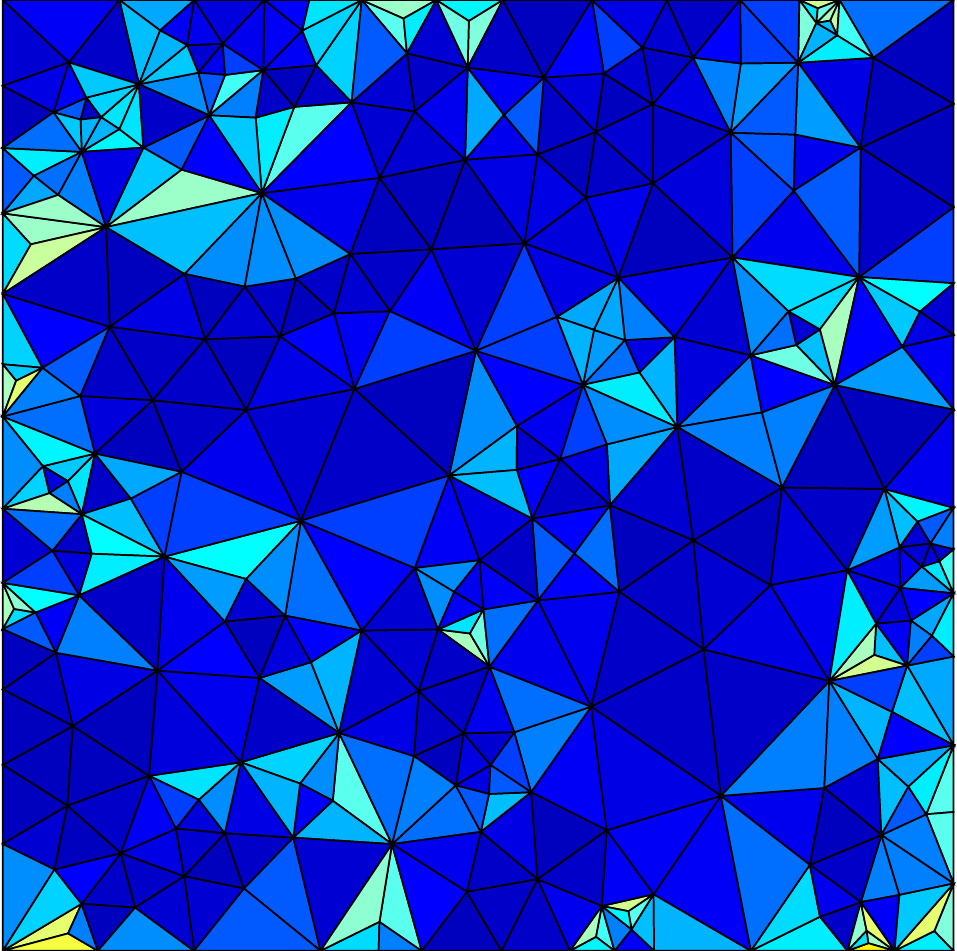
\def\svgwidth{\textwidth/3-1.5mm}
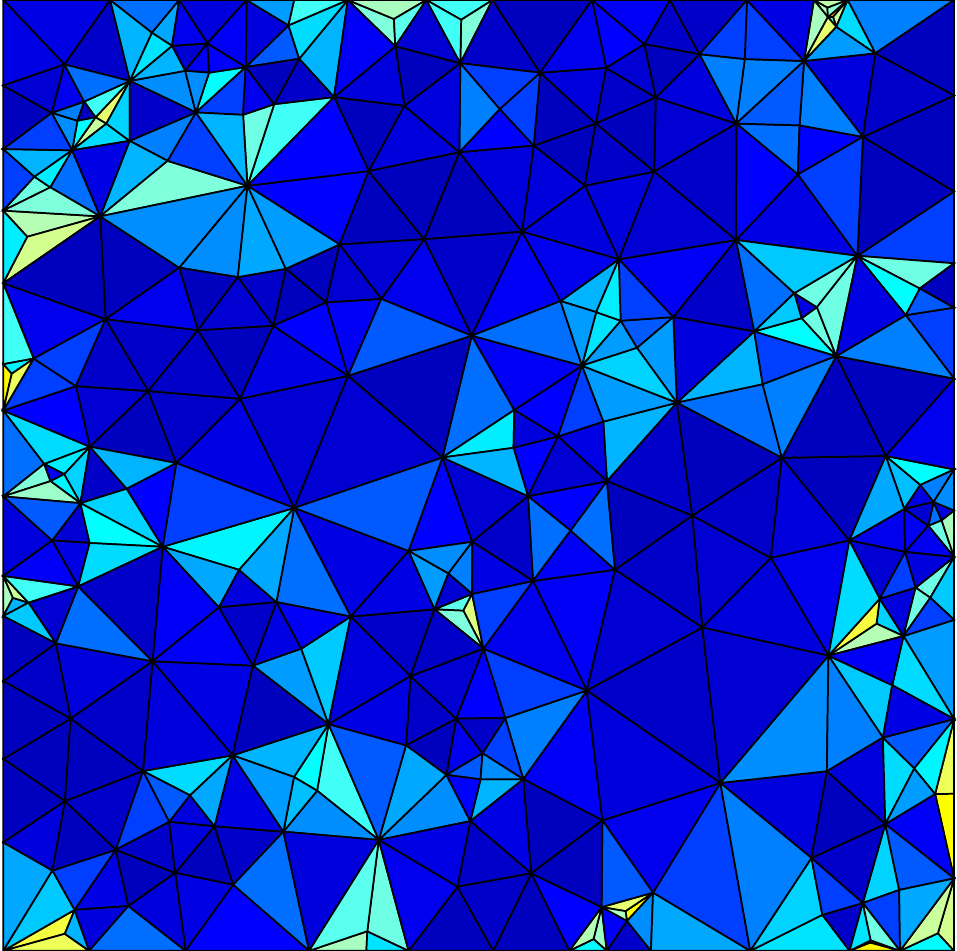\\[1.5ex]
\def\svgwidth{11.4cm}
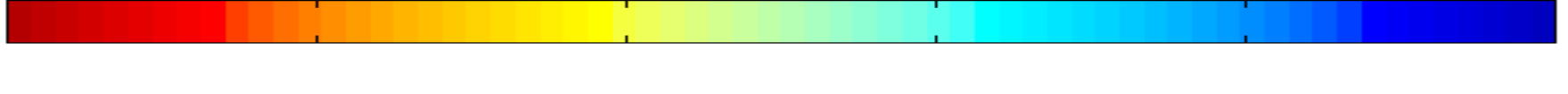\vspace{3ex}
\caption{The initial mesh (left) optimized by $\widetilde\iq_2$ (middle) and by $q_2$ (right). Each element is colored according to its mean ratio quality number.\label{fig_square_exple}}
\end{figure}

In \cite{VartziotisHimpel2013b} we have presented an idea for constructing global optimization-based GETMe smoothing algorithms and hinted in the Outlook at a global optimization-based smoothing method using the isoperimetric quotient. If $E$ is the set of triangles for a planar mesh, $x_e$ are the coordinates of the triangle $e \in E$ and $x_E$ are the coordinates of all nodes in the mesh, then we can define global quality measures in terms of the two-dimensional isoperimetric quotient, for example \begin{equation*}
\iq_2(x_E) = \sum_{e\in E} \iq_2(x_e) \quad \text{and} \quad \widetilde\iq_2(x_E) = \prod_{e\in E} \iq_2(x_e).
\end{equation*}
Let us rewrite
\[
 \iq_2 (x_e) = \frac{1}{C\perim(x_e)^{2}} (\area(x_e) - C\perim(x_e)^{2}) +1,
\]
weigh $\iq_2(x_e)$ by $C\perim(x_e)^{2}$ and shift it, so that all regular triangles attain the same maximal value $0$. This gives a quality measure \[q_2(x_e)=\area(x_e) - C\perim(x_e)^{2},\] which is maximized by regular polygons just like $\iq_2$ is. Moreover, we get a weighted global mesh quality measure
\begin{align*}
q_2(x_E) &= \sum_{e\in E} q_2(x_e).
\end{align*}
While $\iq_2(x_e)-1$ is independent of the size of an element $e$, $q_2$ is more negative for bigger elements of a fixed shape. Therefore the function $q_2(x_E)$ will generally be bigger, if bigger elements are shaped better and smaller elements are shaped worse. Clearly, $q_2$ is concave on the interior nodes and can therefore be used for convex optimization, while giving efficient GETMe untangling and smoothing procedures by way of its gradient ascent. The contribution of $\area(x_E)$ does not depend on the interior nodes of a planar mesh $E$, therefore the gradient of $\area$ is zero on the interior nodes. If $x_{\mathring E}$ corresponds to the collection of variables for the interior nodes $\mathring E$ of $E$, then
\[\nabla q_2(x_{\mathring E}) = -2 \sum_{e\in \mathring E}\perim(x_e)\nabla\perim(x_e).\] 
In Figure~\ref{fig_square_exple} we see the initial mesh, the results from optimizing $\widetilde\iq_2(x_E)$ and from optimizing $q_2$, where the boundary nodes of the square have simply been projected back to the original geometry. We observe, that the optimized meshes look very similar, and that worse triangles tend to be smaller and better triangles tend to be bigger. There is a big difference in usability though. While $\iq_2$ needs a valid and only mildly perturbed mesh, $q_2$ has a unique maximum. It is impressive, but not surprising, that, as long as the mesh topology is the same, $q_2$ will always converge to the same mesh, no matter how the initial coordinates of the interior nodes are chosen.

\begin{figure}[ht]
\centering
\def\svgwidth{\textwidth/3-1.5mm}
\includegraphics[width=\textwidth/3-1.36mm]{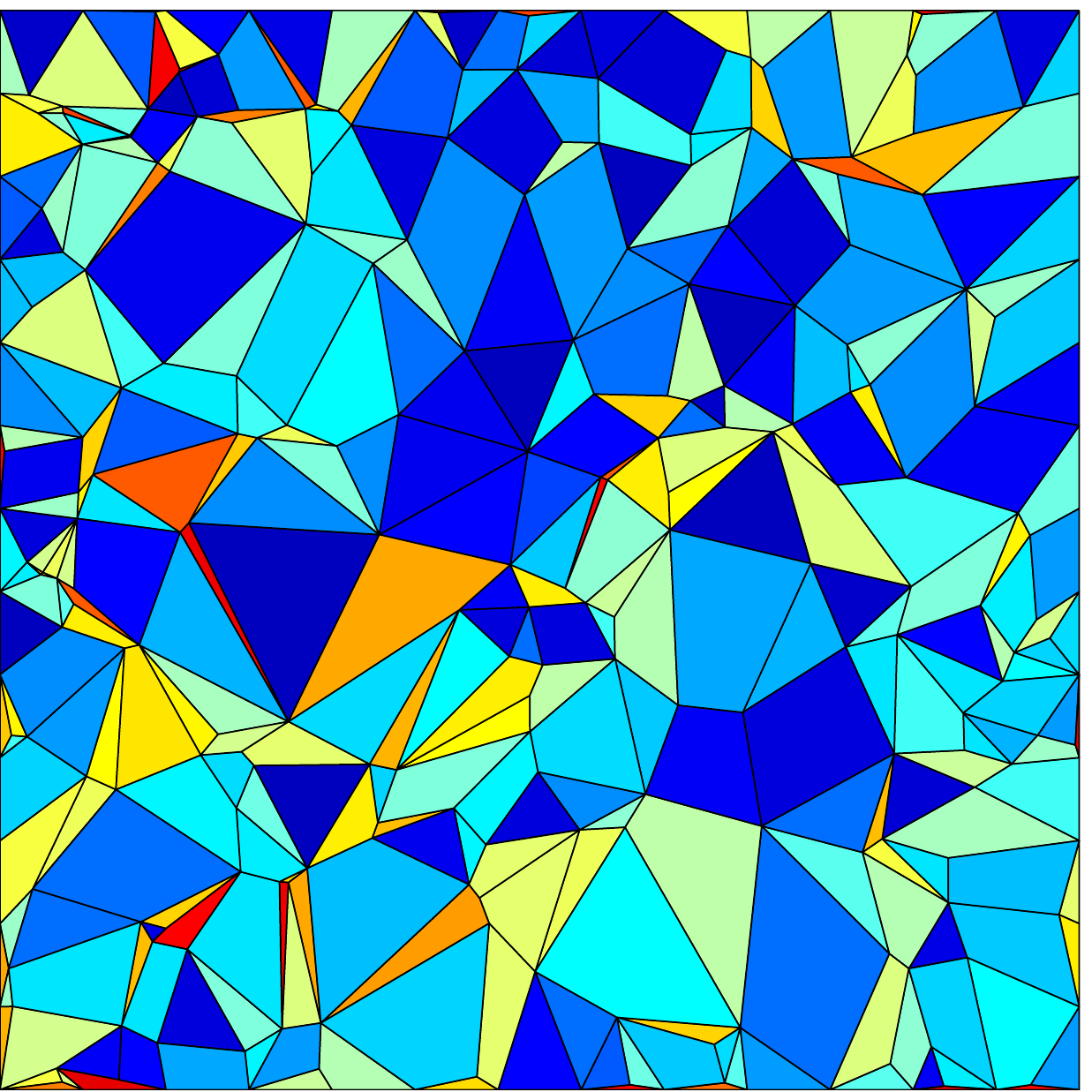}
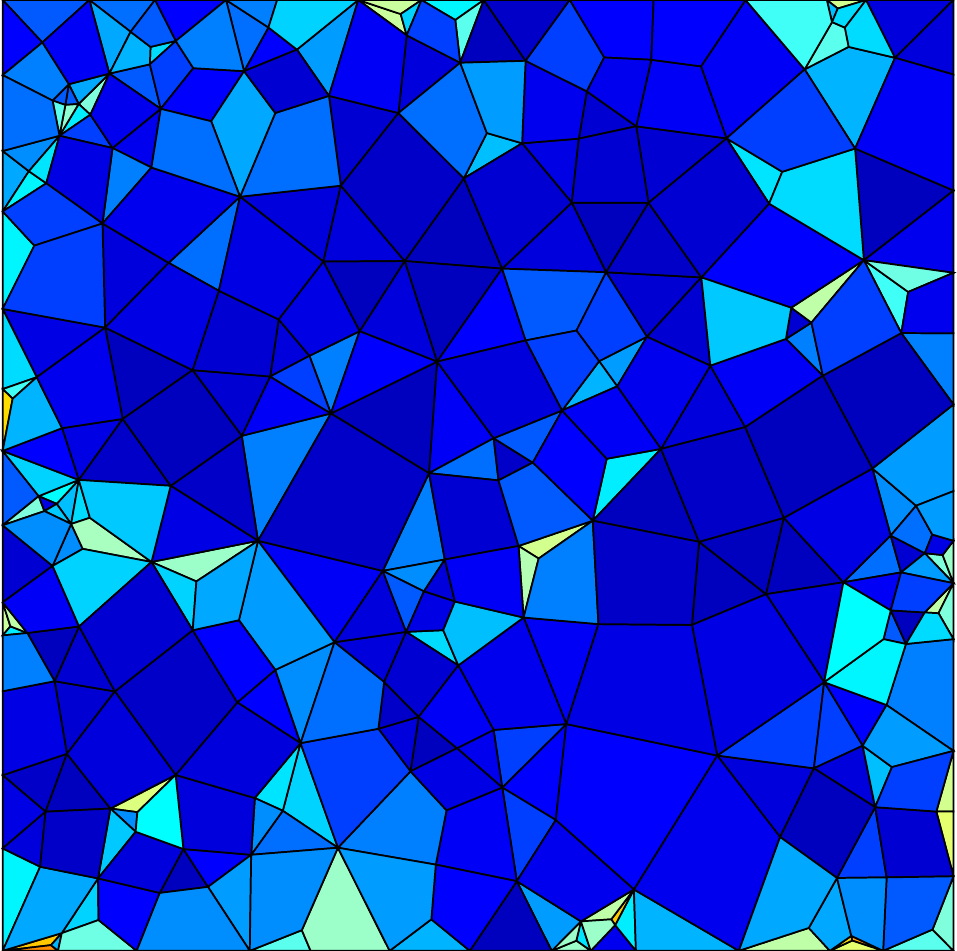
\includegraphics[width=\textwidth/3-1.36mm]{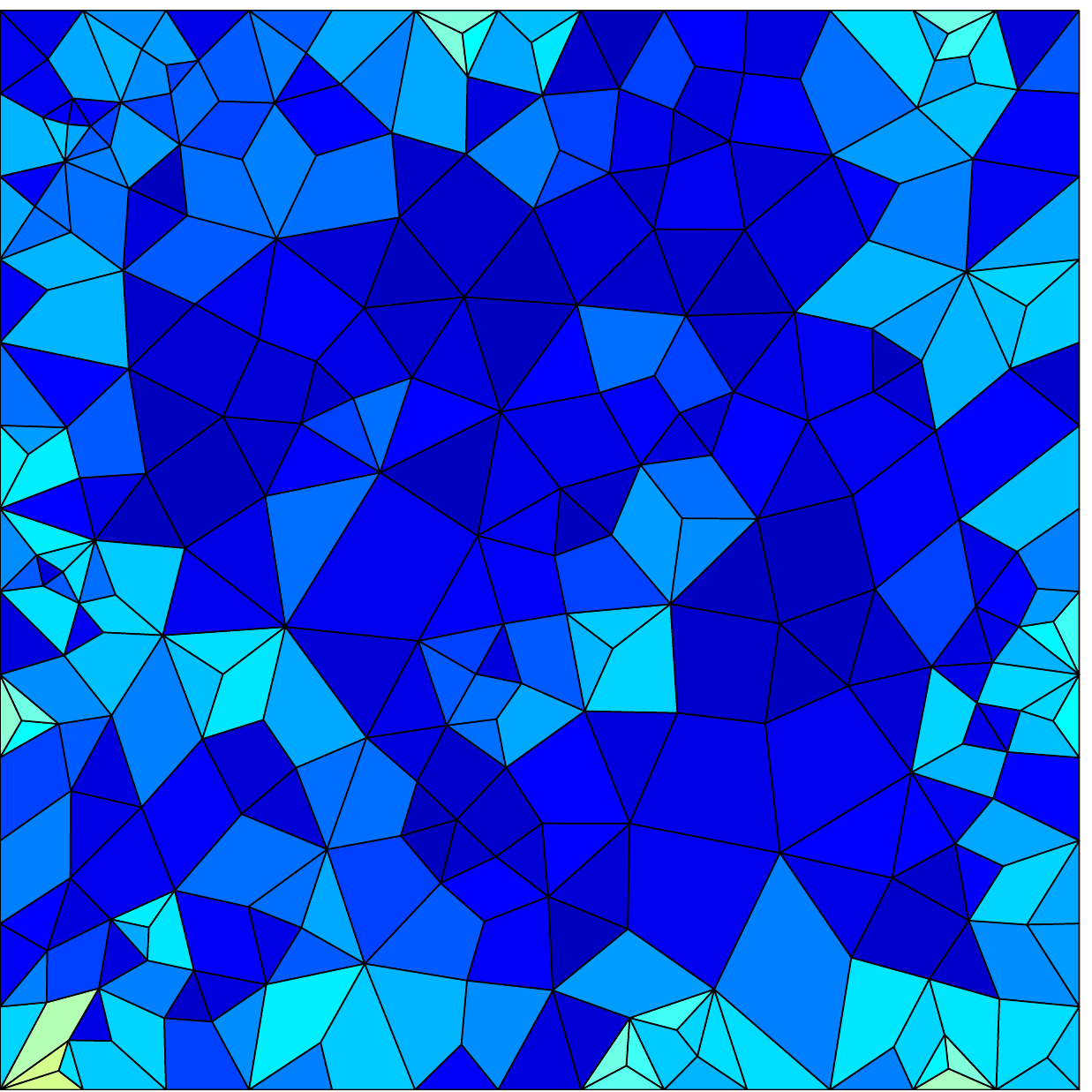}
\caption{The initial mesh (left) optimized by $q_2$ (middle) and by Mesquite (right).\label{fig_square_mixed}}
\end{figure}

The result of smoothing a mixed mesh in comparison to the Mesquite smoothing result can be seen in Figure \ref{fig_square_mixed}. Even though it is not surprising, that the optimization of Mesquite looks better, because it is a global optimization-based method designed to optimize mean ratio, it should be stressed, that the smoothing result by $q_2$ is the unique maximum of the function. The transformation induced by $\nabla q_2$ is able to untangle arbitrary meshes.

\subsection{Smoothing of meshes in three dimensions}

Under the simplifying assumption, that within a given feature area surface meshes in three dimensions should be almost planar, we can use the gradient ascent of $q_2$ to smooth the mesh in combination with the projection to the initial geometry. The gradient of the area can be neglected for almost planar meshes.

If $E$ is a polyhedral mesh, the quality function
\[
q_3(x)=\vol(x) - C^{-1}\area(x)^{3/2} 
\]
and its properties are in exact analogy to $q_2$. If the boundary surface is fixed, $\vol$ is constant and $q_3$ is clearly concave. In practice we might want to smooth the boundary surface by either using a preprocessing step using $q_2$ or smoothing it directly using $q_3$ and projecting it back to the initial geometry. These two functions are only approximately concave, in the sense that the boundary surface will not change much in area and the volume of the three dimensional mesh will not change much as long as the geometry is preserved. 

\begin{figure}[ht]
\centering
\includegraphics[height=2.2cm]{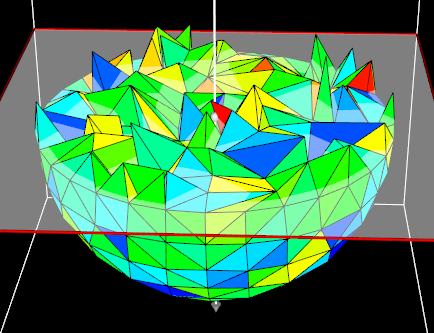}
\includegraphics[height=2.2cm]{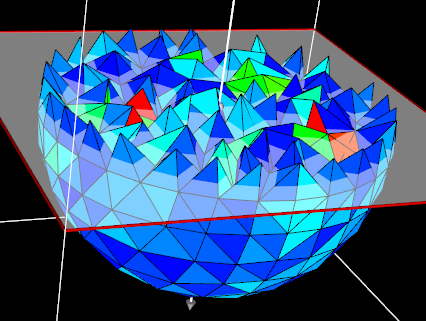}
\includegraphics[height=2.2cm]{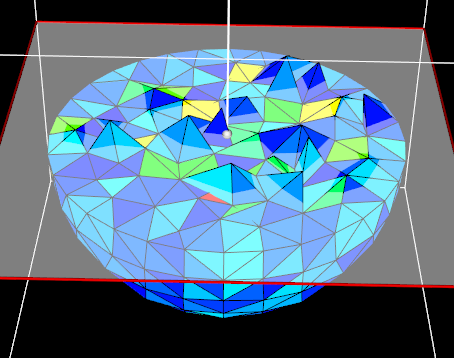}
\includegraphics[height=2.2cm]{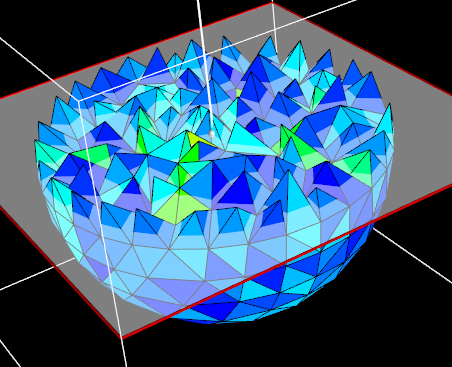}
\caption{From left to right: The initial tetrahedral mesh representing a ball and cut open, the mesh optimized by the original mean ratio, by its square root and by $q_3$.\label{fig_ball_exple}}
\end{figure}

We have tested the smoothing based on $q_3$ on a tetrahedral mesh representing a ball with fixed boundary and compared it with the result of (locally) maximizing the average of the mean ratios via a simple backtracking line search. Even though the average quality attains a very good value in both cases, all of the resulting meshes contain inverted tetraehedra. While the square root of the mean ratio $\sqrt{\mr}$ shows a very good result regarding the mean ratio average, the original mean ratio $\mr$ seems to get stuck in a local optimum. The mean ratio qualities are recorded in Table \ref{table_qualities} and the resulting meshes are depicted in Figure \ref{fig_ball_exple}. The global optimum of $q_2$ is better than the local optimum of $\mr$ and seems to be almost as good as the local optimum of $\sqrt{\mr}$.

Despite the inverted tetrahedra in the resulting meshes, the basic functionality of the gradient ascent of $q_3$ is impressive when compared to the mean ratio. We are guaranteed to find the unique maximum of $q_3$, which also has a good mean ratio average. Just like $q_2$ the three-dimensional analogue $q_3$ is concave, when the boundary is fixed: If $x_{\mathring E}$ corresponds to the collection of variables for the interior nodes of $E$, then
\[\nabla q_3(x_{\mathring E}) = -\sum_{e\in \mathring E}\frac{3}{2C}\sqrt{\area(x_e)}\nabla\area(x_e)^{3/2}.\] In other words, if we keep the boundary nodes fixed in order to preserve the geometry of the mesh, the resulting function is concave, because $\area(x)$ and therefore $\area(x)^{3/2}$ is convex.

In practice, we also want to smooth the boundary mesh. This can be done by using the gradient ascent of $q_2$ or $q_3$ on the boundary and projecting them back to the initial geometry. Even though these functions are not concave on the boundary, they are approximately concave, when the initial geometry is being preserved. If we follow the original GETMe methodology, we would want to make the transformation itself scaling invariant by dividing $\nabla q_3 (x_E)$ or each $\nabla q_3 (x_e)$ by the square root of its norm. In the first case, the transformation would lead to the same mesh.

\section{Laplacian smoothing}

Laplacian smoothing is by far the most popular smoothing method due to its simplicity and time efficiency. Despite its long history, the original Laplacian smoothing has been presented as a heuristic method almost everywhere in the engineering literature, see for example \cite{CanannTristanoStaten1998,Freitag1997,Mukherjee2002,ShontzVavasis2003}.
However, Laplacian smoothing can be derived from a finite difference approximation of the Laplace operator \cite{BuellBush1973}. In particular, it efficiently minimizes a certain convex mesh quality function with a guaranteed and unique result. 
Since we have found very few mentions of it minimizing a simple quadratic energy functional \cite{Ohtake2000,Knupp1999MatrixNorms,Knupp2001}, we will first review the relationship of Laplacian smoothing to the gradient descent of a convex objective function, before we relate it to the popular mean ratio quality criterium and discuss suitable generalizations to polygonal and polyhedral meshes.

\subsection{Laplacian smoothing is global optimization-based}

Let $E$ be the set of triangles and $V$ be the set of vertices of a triangular mesh. Let $V(v) \subset V$ be the set of vertices, which are connected to $v \in V$ by an edge. Then Laplacian smoothing is given by iteratively applying the transformation 
\begin{equation}\label{eq_laplacian_smoothing}
 v \mapsto \frac{1}{|V(v)|} \sum_{v'\in V(v)} v'
\end{equation}
to all inner vertices $v \in V$ of the mesh. Let $x_e$ be the coordinates of the triangle element $e \in E$. 
Consider the convex function
\[
\lambda(x_E) =  \frac{1}{4}\sum_{e\in E} \lambda(x_e), \quad \text{where }\lambda(x_1,x_2,x_3)  = \sum_{i=1}^3 \|x_{i+1}-x_i\|^2,
\]
which can be viewed as a variation of the triangle's perimeter. In fact, $\lambda$ is equal to the square of the Frobenius norm of the (Jacobian) matrix of the difference vectors corresponding to a triangle \cite{Knupp2001}. If we fix the boundary nodes, we can rewrite the above as a sum of edge-length squared over all edges in a mesh, which are not free (i.e. not boundary edges). If we fix all vertices except one inner vertex $v$, then the gradient of $\lambda$ is given by
\[
 \nabla\lambda(v) = |V(v)| v - \sum_{v'\in V(v)} v'.
\]
The critical point of $\lambda$ is the unique minimizer and given by the Laplacian smoothing \eqref{eq_laplacian_smoothing}. Notice, that it is also the gradient descent with scaling parameter $\frac{1}{|V(v)|}$. Therefore, iteratively applying the Laplacian smoothing to all points minimizes the convex function $\lambda$.

We can go further and consider a (modified) isoperimetric quotient \begin{equation}\label{eq_iq}\iq(x_e) = \frac{\area(x_e)}{C\lambda(x_e)} = \frac{1}{C\lambda(x_e)}\left(\area(x_e)-C\lambda(x_e)\right) +1,\end{equation} where $C$ is chosen, such that $\max  \iq(x_e) = 1$ and area is negative on inverted elements. The function $\iq$ is equal to the mean ratio quality criterium in the case of triangles \cite{Knupp2001}.

Notice that the mesh quality function \[\iq(x_E) \coloneqq \sum_{e\in E} \iq(x_e)\] is maximized if and only if \[\sum_{e\in E} (\iq(x_e)-1)\] is maximized. We can view $\iq(x_e)-1$ as a quality measure equivalent to $\iq(x_e)$. Introducing the additional weight $C\lambda(x_e)$ for each element $e$, we get the mesh quality function \begin{equation}\label{eq_objfunction}
q(x_E) \coloneqq (\area-C\lambda)(x_E) = \sum_{e \in E} (\area(x_e)-C\lambda(x_e))
\end{equation}
Clearly the function $\iq(x_e)$ is maximized precisely when $\area(x_e)-C\lambda(x_e)$ reaches its maximal value $0$, so a regular triangle still maximizes the weighted quality function. As long as the boundary nodes are fixed, the $\area$ is constant on a triangular mesh in $\R^2$, so that the problem of maximizing $q(x_E)$
is equivalent to minimizing $\lambda(x_E)$, which can be solved by iteratively applying Laplacian smoothing to all vertices. It even untangles meshes, because the resulting mesh only depends on the mesh topology. This also means, that Laplacian smoothing can result in tangled meshes for certain meshes. In any case, Laplacian smoothing looses its heuristic nature entirely, if we take into account $q$. The result of Laplacian smoothing is clearly related to the mean ratio quality measure weighted by $C\lambda(x_e)$ and shifted so that every regular triangle gives the same maximal value 0.

On a triangular surface mesh in three-dimensional space, vertices moved by Laplacian smoothing will get projected back to the initial geometry and the area of the whole mesh should be almost constant. Therefore, Laplacian smoothing applied in a geometry preserving fashion maximizes the function in \eqref{eq_objfunction}. 

\subsection{Polygonal generalization}

In order to extend Laplacian smoothing to polygonal meshes, we choose
\[
\lambda(x_E) =  \frac{1}{2}\sum_{e\in E} C_e\lambda(x_e), \quad \text{where }\lambda(x_1,\ldots,x_k)  = \sum_{i=1}^k \|x_{i+1}-x_i\|^2\;.
\]
The constant $C_e$ is chosen, so that $\area(x_e)-C_e\lambda(x_e)$ vanishes for a regular polygon $x_e$. For example, $C_e=\sqrt{3}/6$ for a triangular element $e$ and $C_e=1/2$ for a quadrilateral element $e$. Notice, that Laplacian smoothing still maximizes \eqref{eq_objfunction}. If we consider polygonal meshes of mixed type, Laplacian smoothing needs to be adjusted according to the polygon type of the adjacent polygons. Let $E(v,v')$ be the set of (two) elements containing the vertices $v$ and $v'$. \begin{figure}[ht]
\def\svgwidth{2cm}
\centering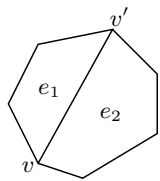\caption{The set $E(v,v')=\{e_1,e_2\}$ of elements containing the vertices $v$ and $v'$.\label{Evw}}
\end{figure} If $v$ is a free inner vertex and all other vertices are fixed, then 
\[
 \nabla\lambda(v) = \sum_{v'\in V(v)}\sum_{e\in E(v,v')}\left(C_e |V(v)| v - C_e v'\right)\;.
\]
Therefore, the correct Laplacian smoothing for mixed surface meshes can be computed by setting the right-hand side to zero. We get the following weighted Laplacian
\[
v\mapsto \frac{1}{|V(v)|}\frac{\sum_{v'\in V(v)}\sum_{e\in E(v,v')} C_e v'}{\sum_{v'\in V(v)}\sum_{e\in E(v,v')} C_e}\;,
\]
which can be rewritten as
\begin{align*}
&v\mapsto \frac{1}{|V(v)|}\frac{\sum_{v'\in V(v)} w_{v,v'} v'}{c_v}\;,\\
\text{where} \quad  &c_v=\sum_{v'\in V(v)}\sum_{v'\in V(v)}w_{v,v'} \text{ and }w_{v,v'} = \sum_{e\in E(v,v')} C_e.
\end{align*}

\begin{figure}[ht]
\centering
\def\svgwidth{\textwidth/3-1.5mm}
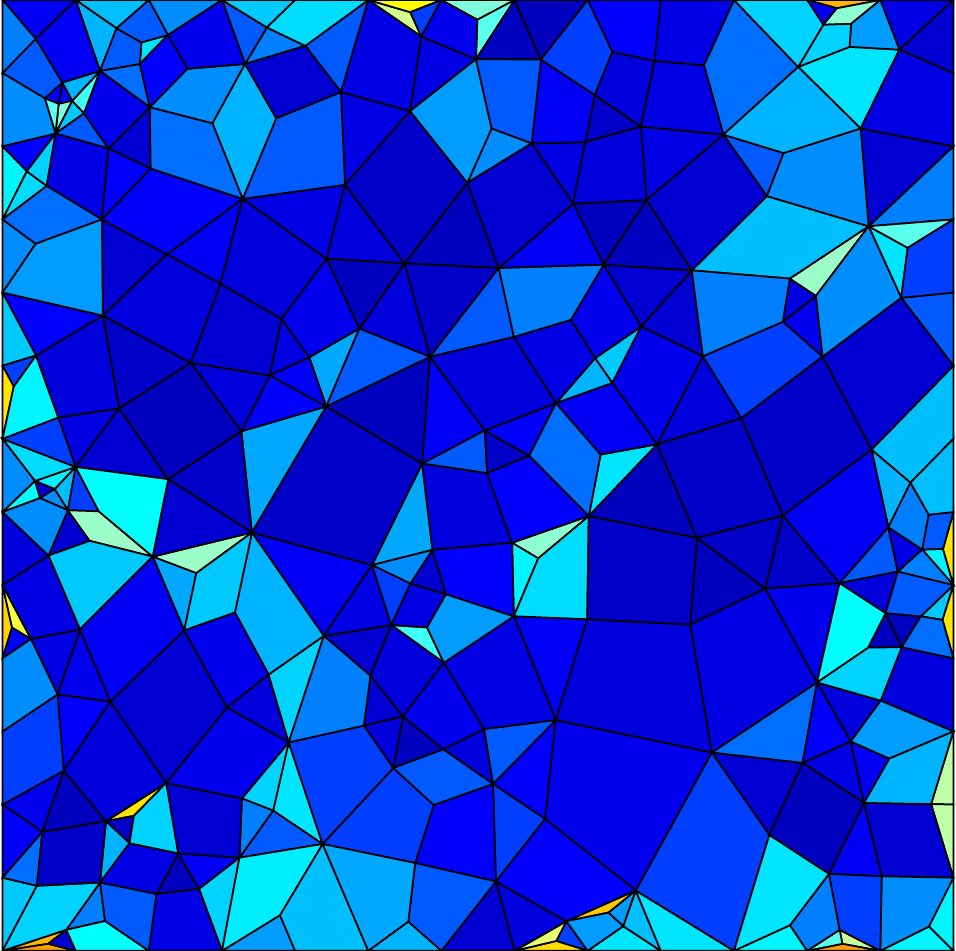
\def\svgwidth{\textwidth/3-1.5mm}
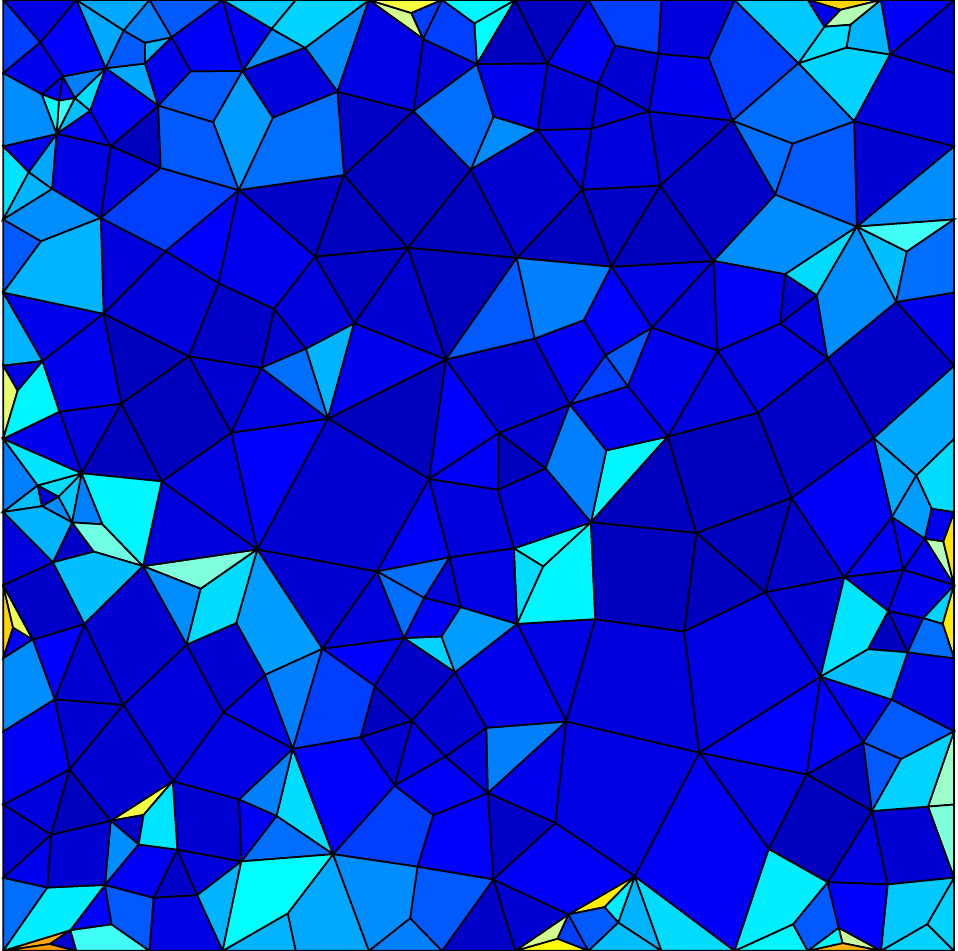
\caption{The result of the unweighted (left) and the weighted (right) Laplacian.\label{fig_square_laplacian}}
\end{figure}

The resulting meshes after applying Laplacian smoothing, weighted and unweighted, can be seen in Figure \ref{fig_square_laplacian}. Clearly, triangles are better in the weighted version.

\subsection{Polyhedral generalization}

Just like in the two-dimensional case, Laplacian smoothing for polyhedral meshes minimizes the sum of edge-length squared
\[
 \lambda(x_E) = \frac{1}{2} \sum_{i<j} \| x_i - x_j \|^2 \quad \text{over all edges in the mesh}.
\]
Just like before the minimizer of $\lambda$ is equal to the maximizer of
\[
 \vol(x_E) - \lambda(x_E),
\]
where we fix the boundary nodes. This cannot be written as a sum of element qualities over all tetrahedra, because the number of attached tetrahedra varies. In retrospect, the original Laplacian smoothing for surface meshes is simply a lucky coincidence, because each inner edge had exactly two adjacent triangles. A better polyhedral generalization to a convex function based on element quality measures is bound to be a bit more involved. 

Instead of starting with the formula for Laplacian smoothing, we can consider a three-dimensional dimensional analogue of \eqref{eq_iq} in order to find a suitable analogue of the convex function \eqref{eq_objfunction}. Again, we would have the decisive advantage that the volume function for a polyhedral mesh is constant, as long as the boundary is fixed.

Firstly, we need to consider an isoperimetric quotient or the mean ratio measure. In addition to the volume function we therefore need to have a homogeneous polynomial of degree three measuring the two-dimensional perimeter. Here we can use any combination of the triangular areas and the edge lengths, hence there are a lot of possibilities, for example
\[
\lambda_1(x_1,x_2,x_3,x_4) = \sum_{i=1}^{4} \area(\hat x_i)  \perim(\hat x_i),
\]
where $\hat x_i$ is the triangle opposite of $x_i$. We then consider the function
\[
 \vol - \frac{1}{C}\lambda_1,
\]
where $C$ is chosen such that \[\max\left(\vol - \frac{1}{C}  \lambda_1\right) = 0.\] The transformation induced by the gradient ascent of the function would then give a result analogous to Laplacian smoothing for surface meshes. The gradient of $\lambda_1$ is determined by its partial derivatives. Instead of $\lambda_1$ we can choose any combination of area and edge length associated to each triangle, which gives a convex function. A non-exhaustive list of examples and their gradients
are listed in Table \ref{table_lambda}.
\begin{table}[ht]
\centering
\begin{tabular}{| >{$\displaystyle}l<{$} | >{$\displaystyle}l<{$} |}
\hline
\centering
\rule{0pt}{2.5ex}f(x)=f(x_1,x_2,x_3,x_4) & \nabla f\\[0.5ex]
\hline
\rule{0pt}{4.5ex}
\lambda_1(x) = \sum_{i=1}^{4} \area(\hat x_i) \perim(\hat x_i) & \sum_{i=1}^{4} \perim(\hat x_i)  \nabla\area(\hat x_i)+ \area(\hat x_i)\nabla\perim(\hat x_i)\\
\rule{0pt}{4.5ex}
 \lambda_2(x) = \sum_{i=1}^{4} \area(\hat x_i)^{3/2} & \sum_{i=1}^{4} \area(\hat x_i)^{1/2}\nabla\area(\hat x_i)\\
\rule{0pt}{4.5ex}
\lambda_3(x) = \left(\sum_{i< j} \| x_i - x_j \|^2\right)^{3/2} & \lambda_3(x)^{1/3} \left(2  \sum_{j\neq k} (x_k-x_j)\right)_{k=1,\ldots,4}\\[0.5ex]
\rule{0pt}{4.5ex}
\lambda_4(x) = \sum_{i< j} \| x_i - x_j \|^3 & \left(3 \sum_{j\neq k} \| x_k - x_j \| (x_k-x_j)\right)_{k=1,\ldots,4}\\
\rule{0pt}{3ex}
 \lambda_5(x) = \area(x)^{3/2} & \area(x)^{1/2}\nabla\area(x)\\
\hline
\end{tabular}
\caption{Different choices of $\lambda$ and its gradient.\label{table_lambda}}
\end{table}

The gradient descents of the $\lambda_i$ result in meshes with excellent mean ratio average, see Table \ref{table_qualities}. Notice, that all of the tested methods resulted in meshes with invalid tetrahedra.

\subsection{Numerical tests}

It comes as no surprise, that the usual Laplacian smoothing given by averaging adjacent vertices for each mesh node is much faster than the gradient ascent methods. The qualities of the resulting meshes speak for themselves. All methods based on $\lambda_i$ give comparable results, which are close to the optimal mean ratio average. It is interesting, that the optimal solution with respect to the mean ratio average is a mesh containing invalid tetrahedra.
\begin{table}[ht]
\centering
\begin{tabular}{| >{$\displaystyle}l<{$} | >{$\displaystyle}l<{$} |>{$\displaystyle}l<{$} |}
\hline
\centering
\rule{0pt}{2ex}\text{Smoothing method} & \text{average } \mr & \text{maximum }\mr \\
\hline
\rule{0pt}{1ex}
\text{Initial} & 0.488795 & 0.986198 \\
\rule{0pt}{2ex}
\mr & 0.706305 & 0.996787\\
\rule{0pt}{2ex}
\sqrt{\mr} &  0.764297 & 0.992954\\
\rule{0pt}{2ex}
\text{Laplacian} & 0.751264 & 0.998345 \\
\rule{0pt}{2ex}
\lambda_1 & 0.762899 & 0.99335 \\
\rule{0pt}{2ex}
q_3,\lambda_2 & 0.761756 & 0.993471 \\
\rule{0pt}{2ex}
\lambda_3 & 0.761317 & 0.994395 \\
\rule{0pt}{2ex}
\lambda_4 & 0.761534 & 0.994798 \\
\rule{0pt}{2ex}
\lambda_5 & 0.761755 & 0.993457 \\
\hline
\end{tabular}
\caption{Mesh qualities with respect to different smoothing methods.\label{table_qualities}}
\end{table}

\section{Enhancing the transformation}

\subsection{Introducing weights}

We can enhance $q_3$ (and similarly $q_2$) as follows. Rather than normalizing $\iq_3$ to equal 1 on regular volume elements, we can let $C$ be arbitrary. It follows that smaller $C$ will make $q_3$ even more negative. Let us therefore impose weights $w_e$ for each element $e \in E$ and consider
\begin{align*}
q'_3(E) &= \sum_{e\in E} q_3(x_e), \quad\text{where } q_3(x) = \vol(x) - w_e\area(x)^{3/2}.
\end{align*}
The new function $q'_3$ is still concave on the inner nodes, moreover $\iq_3(x_e)$ for some $e\in E$ will be bigger in the limit mesh if $w_e$ is increased, because increasing $w_e$ causes $q_3'(x_e)$ to decrease. This way, we gain a control mechanism for the quality of the individual elements in the limit mesh while preserving the convexity of the function. In other words, if we are unhappy with the limit mesh, we can adjust the weights. In particular, if an element $e\in E$ is invalid or of low quality, we can increase $w_e$.

\begin{figure}[ht]
\centering
\def\svgwidth{\textwidth/3-1.5mm}
\input{mesh_iqdiff_normal.pdf_tex}
\def\svgwidth{\textwidth/3-1.5mm}
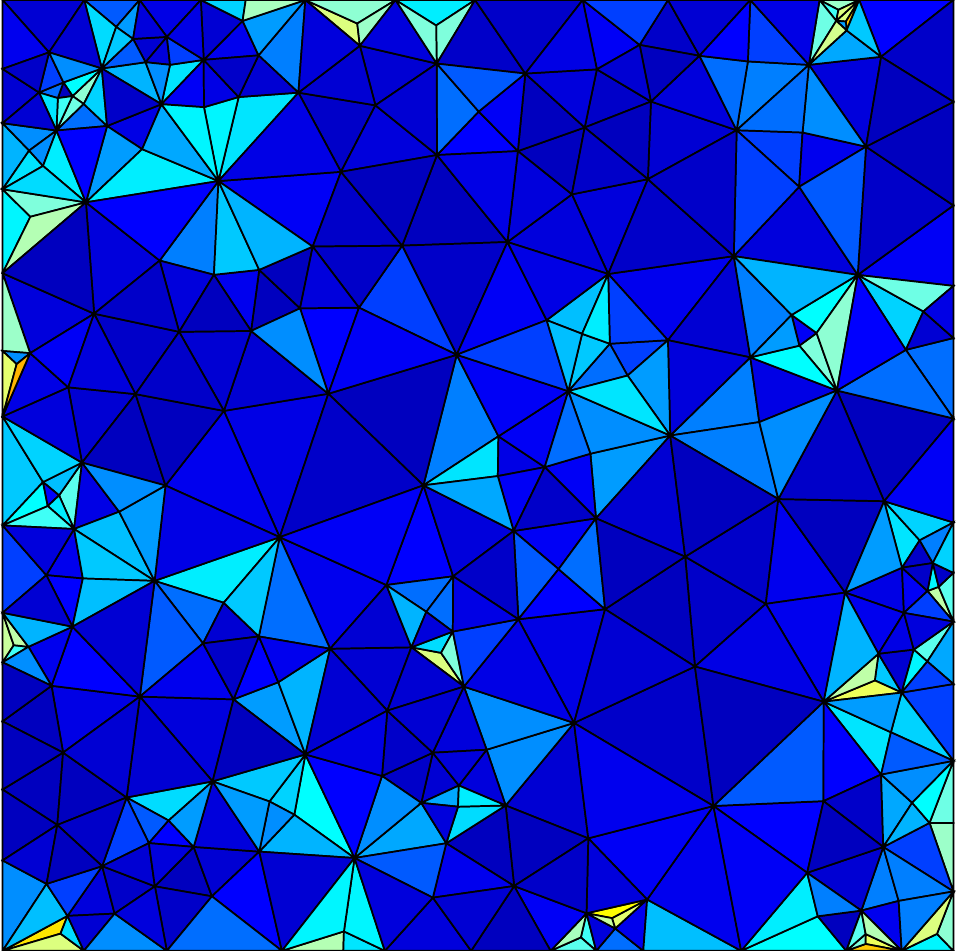
\caption{The mesh optimized by $q_2$ (left), and the mesh optimized by $q_2$ weighted adaptively in order to increase the minimal $\iq_2(e)$ as well as $\iq(E)$.\label{fig_square_diff}}
\end{figure}

Figure~\ref{fig_square_diff} shows the result from optimizing $q_2(E)$, which makes smaller triangles worse and bigger triangles better. 
We can adaptively adjust weights in order to increase the minimal and average value of $\iq_2$. In Figure~\ref{fig_square_diff}, we chose to start with the weight $w_e=1$ for all $e\in E$, adjust the weight as follows:
\begin{enumerate}
 \item Compute $E' = E+ 10^{-3} \cdot \nabla \iq_2(E)$ and let $e'\in E'$ denote the transformed element corresponding to $e \in E$.
 \item If $\iq_2(e')-\iq_2(e)>10^{-5}$ and $\iq_2(e)<0.6$, then set $w_{e'} = 0.99\cdot w_e$.
 \item After adjusting $w_e$ for all $e\in E$, normalize $w_E$ so that $|E|^{-1}\sum w_e =1$.
\end{enumerate}
The above algorithm is heuristic, works okay in this particular example, and can be improved. Nevertheless, for each choice of $w_E$ the weighted function $q'_2$ is convex. We can freely play with the weights in order to optimize the desired combination of minimal and average isoperimetric quotient and the mean ratio.

Weights are usually not necessary, but they could be helpful in particularly bad situation, where the limiting mesh contains invalid elements, like in our example of the ball represented by a tetrahedral mesh. However, adaptively changing weights might effectively mean, that we are not doing convex optimization anymore. Topology modifications might be a more practicable solution.

\subsection{Topology modifications}

Due to the scaling dependence of $q$, bad elements tend to be smaller. Therefore, topology modification is more than ever a necessary evil. However, smaller elements can be removed more easily. We provide tests, that show the effectiveness of edge-swap, edge-collapse and vertex-split in combination with the new convex optimization-based smoothing method.

\begin{figure}[ht]
\centering
\def\svgwidth{\textwidth/3-1.5mm}
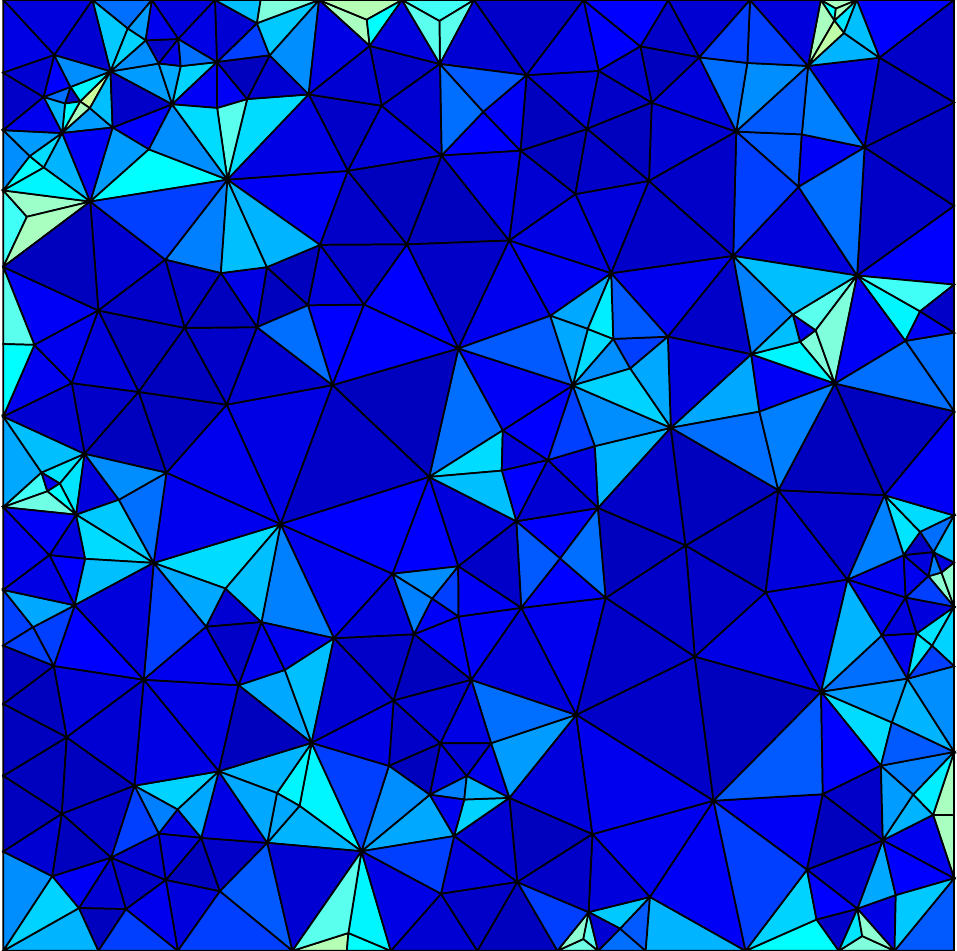
\def\svgwidth{\textwidth/3-1.5mm}
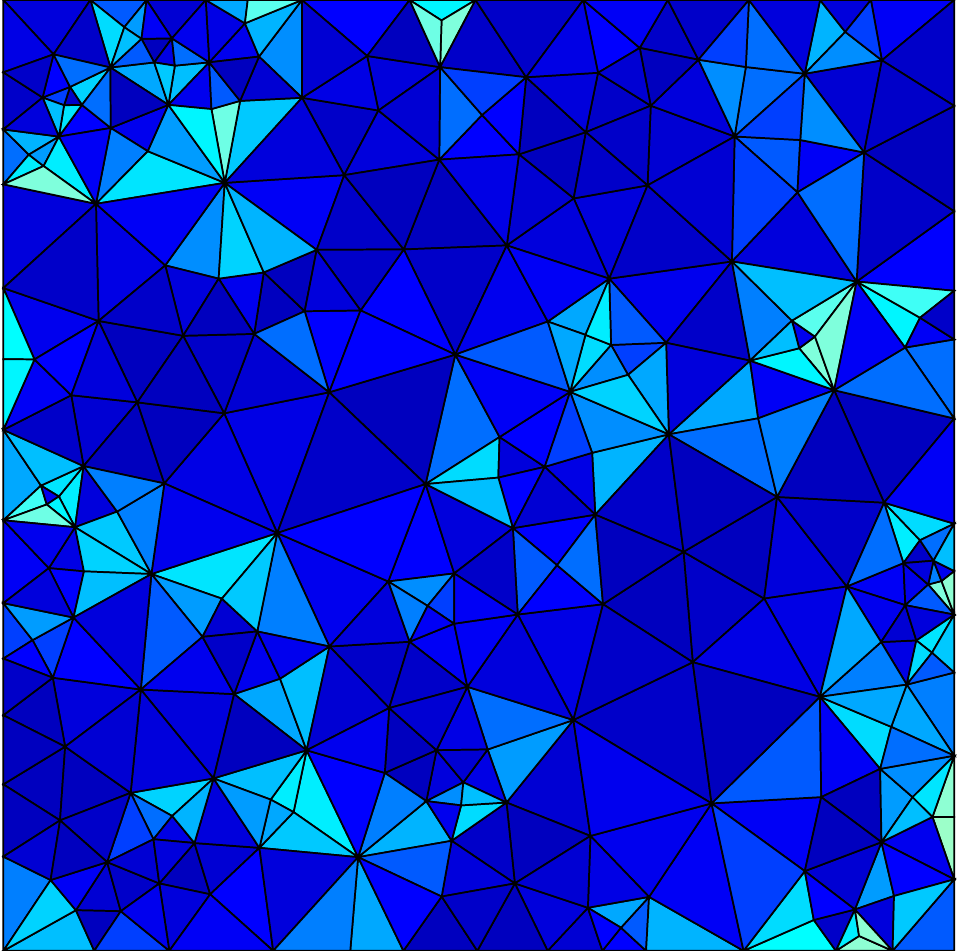
\caption{The intermediate and final modification result deleting all triangles with $\iq_2(e) < 0.6$.\label{fig_square_mod}}
\end{figure}

We have seen, that minor adjustments to $q'_2$ can be made in order to improve the quality measured by $\iq_2$. However, in order to use this approach to industrial problems, we can try to utilize the special behavior of $q_2$ with respect to the perimeter of the triangles, namely that it gives preference to bigger triangles. More precisely, when triangles become are badly shaped, we should modify the connectivity in order to delete them. Unless the corners of the triangles are constrained to some feature edge, these triangles tend to be small and easily removable after optimizing $q_2$. In Figure \ref{fig_square_mod} we modified the connectivity after an initial optimization run (resulting in the left figure of Figure 2) and removed all triangles $e$ with $\iq_2(e)<0.6$ by edge-collapse in addition to optimizing the modified the result. Figure \ref{fig_square_mod} shows an intermediate step and the final result containing only triangles with $\iq_2(e)\ge 0.6$.
\section{Summary and outlook}

We have shown that the driving force behind the Laplacian smoothing is the popular mean ratio quality measure. We have used these ideas to find convex quality functions, whose gradient ascent yield effective untangling and smoothing methods. Such methods have a tendency of resulting in meshes with worse triangles being slightly smaller, which lends itself to topology modifications.

While most smoothing methods try to preserve validity in each smoothing step, we propose to incorporate topology modifications in the smoothing process to remove any invalid or bad elements. We will discuss these ideas in a future paper.

The recent developments on viewing GETMe methods as the gradient ascent of the volume function \cite{VartziotisHimpel2014} and the present work raises the following philosophical question: Does every smoothing method ultimately get its smoothing power from some underlying quality function?
\bibliographystyle{siam}
\bibliography{literature.bib}

\end{document}